\documentclass[11pt]{article}
\usepackage[a4paper]{geometry}

\usepackage[utf8]{inputenc}
\usepackage[T1]{fontenc}
\usepackage{lmodern}
\usepackage[sectionbib,numbers]{natbib}

\usepackage{amstext,amsthm,amsmath,amssymb,mathtools,bbm,mathrsfs}
\usepackage[colorlinks=false]{hyperref}
\usepackage{graphicx}
\usepackage{enumerate}

\usepackage{graphicx}
\usepackage[normalem]{ulem}

\newtheorem{proposition}{Proposition}
\newtheorem{corollary}[proposition]{Corollary}
\newtheorem{theorem}{Theorem}

\theoremstyle{definition}

\newtheorem{example}[proposition]{Example}
\usepackage{color}
\usepackage[colorinlistoftodos,color=lightgray]{todonotes}

\usepackage{txfonts, upgreek}

\newtheorem{xx}{\bf xxx}

\newtheorem{zz}{\bf zzz}

\newtheorem{yy}{\bf yyy}

\definecolor{groen}{rgb}{0,0.5,0.2}

\newcommand{\ts}{\hspace{0.5pt}}

\newcommand{\II}{\mathbb{I}}

\newcommand{\KK}{\mathbb{K}}

\newcommand{\NN}{\mathbb{N}}
\newcommand{\Nnull}{\mathbb{N}_{0}^{}}
\newcommand{\PP}{\mathbb{P}\ts}

\newcommand{\RR}{\mathbb{R}}

\newcommand{\ZZ}{\mathbb{Z}\ts}
\newcommand{\Z}{\mathbb{Z}\ts}

\def\cB{\mathcal{B}}

\def\cE{\mathcal{E}}
\def\cF{\mathcal{F}}

\def\cM{\mathcal{M}}

\def\cL{\mathcal{L}}

\def\cP{\mathcal{P}}

\def\cU{\mathcal{U}}

\def\mfU{\mathfrak{U}}
\def\mfu{\mathfrak{u}}

\DeclareMathOperator{\Pois}{Poiss}
\DeclareMathOperator{\Exp}{Exp}

\DeclareMathAlphabet{\mathpzc}{OT1}{pzc}{m}{it}

\newcommand{\GG}{\mathbb{G}}
\newcommand{\UU}{\mathbb{U}}
\newcommand{\wt}{\widetilde}
\newcommand{\ui}{\underline{i}}
\newcommand{\uj}{\underline{j}}
\def\wh{\widehat}
\def\uur{\underline{\underline{r}}}
\def\nto{\mathop{\Longrightarrow}\limits_{n \to \infty}}
\def\Nto{\mathop{\Longrightarrow}\limits_{N \to \infty}}
\def\jto{\mathop{\Longrightarrow}\limits_{j \to \infty}}
\newcommand{\noi}{\noindent}

\begin{document}
\title{\LARGE Stochastic Evolution of spatial populations:\\
From configurations to genealogies and back}
\author{\sc Andreas Greven$^1$}
\date{{\today}\\
}

\maketitle

\begin{abstract}
The paper reviews the results obtained for spatial population models 
and the evolution of the genealogies of these populations during the last decade by the author and his coworkers.
The focus is on their large scale behaviour and on the analysis of universality classes of large scale behaviour via the methods of the hierarchical mean-field limit or via the spatial continuum limit and from another angel, the finite system scheme .
We use genealogical information to analyze the type and location structure and vice versa.

To apply this approach and to explain effects in biological situations we extend the classical model classes in new directions.
Namely we look as population models here at: 
Fleming-Viot genealogies in continuum geographic space, Cannings models with block resampling (reducing diversity), Fisher-Wright diffusion with coloured seedbanks (enhancing diversity) and evolving genealogies of Fleming-Viot models with selection and rare mutation.

\bigskip

\paragraph{Keywords:} 
Spatial population models, genealogies, ultrametric measure spaces, universality, renormalization, finite system scheme, continuum space limit.

\bigskip

\noi {\bf AMS Classification: Primary: 60B10, 05C80; Secondary: 60B05, 60G09}

\noi {\bf MSC: Primary 60J70, 60K35; Secondary 92D25}

\end{abstract}

\bigskip \footnoterule

\noindent {\footnotesize $^{1)}$ Department Mathematik, Friedrich-Alexander Universit\"at
  Erlangen-N\"urnberg, Cauerstr. 11, D-91058 Erlangen, Germany,\\
  greven@math.fau.de}

\newpage

\tableofcontents

\newpage

\section[Introduction]{Introduction: Basic ideas and the models}\label{AG-s.motback}

\textbf{Basic ideas and goals} \quad
We want to present here four main ideas to analyze the behaviour of stochastically evolving spatial populations.
The latter means we have individuals in a \emph{geographic space} $\GG$, for example $\ZZ^d$ or $\RR^d$ and they have genetic \emph{types} from a set $\KK$ and potentially they also have  \emph{internal states} from a set $\II$ which allows for example to model that we have a seedbank where we distinguish active individuals from dormant ones (the ones in the seedbank).
The state of the population is coded via a counting measure on $\GG \times \KK \times \II$.

In the case of individual based models, the population evolves then by individuals \emph{migrating} in geographic space, \emph{mutating} in type space, \emph{switching} the internal state, \emph{change of generation}, i.e. individuals are replaced by their descendants via neutral \emph{resampling}, resampling with \emph{selection} and/or \emph{recombination}, where the type of a descendant is chosen according to some rule we will explain below.

All transitions are stochastic.
Basic references for these classes of models are \cite{D93,EK86}.

We will primarily be interested in frequencies of types at the locations in the limit of large number of individuals (per site) of small mass and suitably scaled parameters leading to a \emph{stochastic} limit evolution, which will be \emph{stochastic} because of \emph{neutral resampling}, all other terms becoming deterministic drift terms.

\paragraph{Outline}
We outline the \emph{goals} of our investigations and the \emph{key models} in this section, in Section 1.2 and 1.3 we describe the \emph{principal mathematical tools} we developed for our work.
In Section 1.4 we explain a number of \emph{typical results} we were able to prove and in Section 1.5 sketch \emph{challenges for the future}.

\subsection{Goals}\label{ss.goals}
Next we formulate the \emph{four main goals} of our work related to the analysis of the large time-space scale behaviour of spatial populations under stochastic evolution.

\medskip

\noindent \textbf{\emph{(i) From measure valued to genealogy valued processes.}} \quad
The first point is now that, if we are interested in the evolution of the frequencies of the types at the different geographic sites, then for the analysis of this \emph{measure-valued} process on $\GG \times \KK \times \II$ it is useful to have \emph{genealogical information} about the population as for example information on the \emph{family tree} of the population alive at the \emph{current time $t$}, more precisely the marked (with type and location) current family tree.
This allows to understand for example the formation of geographic monotype patches forming in low dimensional spatial Fleming-Viot models as time increases, which is a feature of the configuration process.

On the other hand we are in applications often interested a priori in the process of genealogies (law of the genealogical distances of typical individuals, time back to the most recent common ancestor (MRCA), the structure of subfamily decompositions).
Here results one obtains on the type-space occupation measure allow in turn to prove properties of the genealogies by studying marked genealogies and their evolution, see \cite{MAX_Griesshammer}.

For both these reasons we want to characterize the processes of the \emph{evolving} genealogies as solutions of \emph{wellposed martingale problems} whose projection on type-locations i.e. a measure-valued process, are themselves solutions of wellposed martingale problems. 

This requires first that we specify \emph{Polish spaces} which we can use as \emph{state spaces} for these evolving genealogies, specify suitable \emph{test functions} on these spaces and define the \emph{operators} for the martingale problems on these test functions.
Then we need to show that this is indeed leading to a \emph{wellposed martingale problem} and hence is specifying our model.

Here we can use the approach of Depperschmidt, Greven, Winter and Pfaffelhuber based on the use of \emph{equivalence classes of marked ultrametric measure spaces} to describe genealogies  \cite{GPW09,DGP11,GPWmp13,ggr_tvF14}.
For a survey of this approach we refer also to \cite{DG18evolution}.
With this method we started with Depperschmidt, Pfaffelhuber and Sun  analyzing the genealogies in various models of the above mentioned type, which we will describe in the sequel, \cite{DGP13,GdHKK,GSW,GKW}, see also \cite{infdiv,ggr_GeneralBranching} for branching models.
\medskip

\noindent \textbf{\emph{(ii) Universality.}} \quad
The second point is, to use for the analysis of the \emph{universal properties} of the population in \emph{large time-space scales}, that we can look at the evolution mechanisms in multiple time-space(-rate) rescalings  and consider the scaling limit of these \emph{renormalized systems} for \emph{large time-space scales}.

\textit{(ii-a)}\quad 
It now turns out that this analysis is typically too complex.
We therefore choose different geographic spaces namely the \emph{hierarchical groups $\Omega_N$} with $N \in \{2,3,\cdots\}$ instead of $\ZZ^d$ and in particular the interesting $\ZZ^2$.
However we have comparable \emph{potential theoretic properties} for random walks on this geographic space, such that this gives as $N \to \infty$ an approximation of the large scale behaviour on $\ZZ^d$, a procedure called \emph{hierarchical mean-field limit} which was developed by Dawson and Greven in \cite{DG93,DGV95,DG96,CDG04,DGdHSS}.
This way one is then able to find the \emph{universality classes} for the large scale behaviour. See for the treatment of a new class of models explained later \cite{GdHKK,GHK17} and current work on models with seedbank \cite{GdHOpr3}. 
Note however also that in \emph{many} ecological systems $\Omega_N$ is in fact even better suited than $\ZZ^d$ or $\RR^d$.

The problem we focused on in our research was to carry out this approach now on a \emph{genealogical level}, but also extending the analysis of the measure-valued case, to \emph{new} classes of models, for example passing from Fisher-Wright models to the \emph{Cannings model} with a geographically structured Cannings mechanism or allowing in Fisher-Wright models \emph{seedbanks} resp. the classical version but with  selection and recombination  \cite{ghs1,ghs_renorm,GdHKK,GHK17, DGsel14}.

\textit{(ii-b)}\quad 
An alternative method to characterize universality classes of the large time-space scale behaviour of systems on $\ZZ^d$ or $\Omega_N$ is the \emph{continuum space limit} resulting in models on $\RR^d$ respectively the continuum space hierarchical group $\Omega^\infty_N$, whose properties on large time-space scales are typical for all dynamics in the domain of attraction, see \cite{GSW,ghs1}.

\textit{(ii-c)}\quad 
A related question here is the question what aspects of the behaviour of populations in \emph{infinite} geographic spaces are relevant for \emph{large but finite} systems of the real world and how can this be exhibited.
Here the technique is to consider the space-time scales in which, the size of space tends to infinity and time grows as function of the size such that nevertheless the behaviour of the finite system starts displaying the finiteness of the space. 

Indeed the \emph{finite system scheme} has been introduced by Cox and Greven to exhibit what relevance properties of the longterm behaviour of systems on infinite geographic spaces have for large but finite geographic systems.
This has been introduced first for infinite particle systems \cite{CG90} and then interacting diffusions \cite{CGSh95} or measure valued processes \cite{DGV95}.

This is now needed for some \emph{new} systems relevant in biology as Cannings models with blockresampling or Fisher-Wright diffusions with seedbanks and is needed on the level of the \emph{evolving genealogies}, which requires new tools and concepts \cite{GKW,GdHOpr2} and has been considered by for the Fleming-Viot model and the Cannings model.
Another direction is to drop the "constant population size property" and to consider \emph{branching models} and their evolving genealogies, \cite{ggr_tvF14,infdiv,ggr_GeneralBranching}.
\medskip

\noindent \textbf{\emph{(iii) High and low volatility, continuum geographic space.}} \quad
As a third point of our research we aimed at extending the class of spatial population models as Fleming-Viot or branching models to be able to treat a number of effects which are relevant in biological applications. 

First look at change of generation involving \emph{high volatility} due to catastrophic events, leading us to the \emph{Cannings model with block-resampling} resulting in a larger tendency to local monotype blocks.
This leads to situations, where classical models show coexistence of types but block resampling induces monotype clusters, see \cite{GdHKK,GHK17}.
On the other end the opposite effect of extreme \emph{reductions of volatility} and hence increasing genetic diversity is leading us to consider \emph{seedbanks} and to show that we can different from classical models have coexistence of types in $d=2$.

A second point is to be able to work with \emph{continuum} geographic space instead of discrete geographic spaces, which is often more natural.
For genealogies see \cite{GSW} for a discussion and further references.
\medskip

\noindent \textbf{\emph{(iv) Nonlinear effects.}} \quad
The fourth point is that we can apply the points (i) to (iii) to study the evolution of populations through the interplay between rare mutation, subsequent \emph{selection} and the influence of migration from this genealogical perspective, see here \cite{DGsel14,GPPW}.
This way we obtain insight in \emph{spatial sweeps}.
To investigate the associated \emph{genealogical} structure and including \emph{recombination} is however still work in progress.
Similarly this is true for populations under competition for scarce  resources modeled with \emph{logistic branching} models.
\medskip

\subsection{Models}\label{ss.models}
We investigate these points in a class of stochastic evolutions of the following type.
We assume that we have metric spaces $(\GG,r_\GG),(\KK,r_\KK),(\II,r_\II)$:
\begin{align}\label{AG-e179}
& \GG \quad \mbox{is a countably infinite abelian group}, & \\
& \KK \quad \mbox{ is a compact Polish space},  & \\
& \II \quad \mbox{is a finite (or countable) set.} &
\end{align}
In particular we can consider $\GG \times \KK \times \II$ with the product topology and its Borel-$\sigma$-algebra and construct the corresponding \emph{measure valued processes} on this space, generically denoted by $X=(X(t))_{t \geq 0}$.
\medskip

Given are the following ingredients for the dynamics.
Consider for the moment a population with a finite number of individuals per site:
\begin{itemize}
\item \emph{Migration} rates $a(\cdot,\cdot)$ on $\GG \times \GG$, which are homogeneous, i.e. \; $a(i,j)=a(0,j-i) \; ; \; i,j \in \GG$.
This describes the rates at which an individual migrates from $i$ to $j$.
\item A \emph{resampling rate} $d$ and a \emph{resampling measure} $\Lambda$ on $[0,1]$, which is a probability measure regulating the natural change of generation, for $\Lambda=\delta_0$ we talk just about \emph{resampling} and if $\Lambda(\{0\})=0$, we call this \emph{$\Lambda$-resampling}.

In a \emph{resampling event} two individuals are replaced independently by choosing a parent of the two children with probability $\frac{1}{2}$ for each of the two.

In a \emph{$\Lambda$-resampling event} we choose $r$ with intensity $\Lambda^\ast(dr)=r^{-2} \Lambda(dr)$, mark each individual with a $1$ with probability $r$.
Choose one of the $1-$marked individuals at random, his offspring replaces all other $1$-marked individuals.
\item \emph{Mutation} rates $m(\cdot,\cdot)$, a transition kernel from $\KK$ to $\cB (\KK)$, $m(u,dv)$ describes the rate to mutate from type $u$ to a type in $dv$.
\item A \emph{fitness function} $\psi:\KK \times \KK \to [0,1]; inf \; \psi=0, \sup \psi=1$ and $s \geq 0$ the \emph{selection rate}, here $\psi (u)-\psi(v)$ describes the \emph{bias} with which in a resampling event "$u$ is chosen over $v$ as descendant" compared to the neutral case.
\item A \emph{recombination} rate $r$ and a recombination kernel $r(\cdot,\cdot)$ from $\KK \times \KK$ to $\cB(\KK)$.
\item Rates $e$ and $eK$ at which individuals may become \emph{dormant} (i.e. enter a \emph{seedbank}) respectively become \emph{active} again.
\end{itemize}
We shall shortly also discuss branching dynamics (\cite{infdiv,ggr_GeneralBranching,ggr_tvF14}), which allows to have variable population size.
Consider now the Markov jump process $X$ of measures on $\GG \times \KK \times \II$ which specify the number of individuals at a geographic site of a specific type and internal state.
In all these models above we pass to the \emph{type frequencies} per site and then take the \emph{large population limit}.
This results then in a \emph{measure-valued diffusion} or measure valued c\'{a}dl\`{a}g process, which is the model we want to study and characterize by wellposed martingale problems. 
For these models it is our goal to construct in Section~\ref{AG-s.descgen} via wellposed martingale problems the \emph{marked-genealogy process}.

For above processes it is often natural to consider them in \emph{random  medium}, i.e. one of the specified parameters is \emph{spatially fluctuating}.
Then it is of interest how such systems compare to the homogeneous version specified above, see \cite{GHK17}.

The \emph{state spaces} of the processes of frequencies is
\begin{equation}\label{AG-e2100}
S=\left(\cP(\KK \times \II)\right)^\GG \quad,
\end{equation}
where $\cP (E)$ is the set of probability measure on a Polish space $E$. 
This leads to the class of interacting Fleming-Viot processes ($\Lambda$-part not present) respectively interacting $\Lambda$-Cannings models if $\Lambda(\{0\})=0$.
The martingale problem for this dynamics looks as follows.

Let $\cF$ be functions (polynomials) of the form 
\begin{align}\label{AG-e217}
F(\underline{x})= & \int\limits_\KK \cdots \int\limits_\KK f(u_1,u_2,\cdots, u_n) x_{\xi_1}(du_1) \cdots x_{\xi_n} \xi_n(du_n), & \\
\mbox{with  } & f \in C_b(\KK,\RR), \quad \xi_1,\cdots,\xi_n \in \GG, \quad  \underline{x}=(x_\xi)_{\xi \in \GG}, x_\xi \in \cP(\KK). & \nonumber
\end{align}
For each of the mechanisms listed above and acting \emph{independently} we have a corresponding operator, with the obvious names.
Hence we set
\begin{align}\label{AG-e224}
G= G^{\rm mig}+G^{\rm res}+G^{\rm Cann}+G^{\rm mut}+G^{\rm sel}.
\end{align}
Here $G^{\rm mig}$ is the \emph{interaction} between geographic sites, $G^{\rm mut}+G^{\rm sel}$ is the sum over the \emph{first order} terms  of the evolution of the components over the different sites and $G^{\rm res}+G^{\rm Cann}$ are the sums over the \emph{second order} respectively the \emph{integral operator terms} over the different sites, the latter corresponding to the jump process part.
Here are the actions of the single operators.

Let $x=(x_\xi)_{\xi \in \GG} \in (\cP(\KK))^\GG$, then:
For the resampling rate $d$ the resampling operator reads:
\begin{equation}\label{AG-e240}
G^{\rm res}(F)= \sum\limits_{\xi \in \GG} d \cdot \int\limits_\KK \int\limits_\KK \frac{\partial^2 F(x)}{\partial x_\xi \partial x_\xi} (u,v) \; Q_{x_\xi} \; (du,dv),
\end{equation}
with $Q_x(du,dv)=x(du) \delta_u (dv)-x(du)x(dv), x \in \cP(\KK)$.

For $\Lambda^\ast$ with $\Lambda(\{0\})=0$ we define
\begin{align}\label{AG-e262}
&G^{\rm Cann}(F)= \sum\limits_{\xi \in \GG} \sum\limits_{i=1}^n \int\limits_{(0,1]} \Lambda^\ast(dr) & \\
& \qquad \int\limits_\KK x^i_\xi (du) \left(F(x^1_\xi, \cdots,x^{i-1}_\xi, (1-r)x^i_\xi + r \delta_u,x^{i+1}_\xi, \cdots,x^n_\xi)-F(x_\xi)\right).&\nonumber
\end{align}
For this latter transition a spatially structured version will be considered later on.

For mutation rate $m$ and kernel $M$:
\begin{equation}\label{AG-e253}
G^{\rm mut}(F)= \sum\limits_{\xi \in \GG} m \int_\KK \left\lbrace \int_\KK   \frac{\partial F(x)}{\partial x_\xi} (u) M(u,dv)-\frac{\partial F(x)}{\partial x_\xi} (u) \right\rbrace x_\xi(du).
\end{equation}
For selection rate $s$ and selection function $\psi$:
\begin{equation}\label{AG-e257}
G^{\rm sel}(F)= \sum\limits_{\xi \in \GG} s \cdot \int\limits_\KK 
\frac{\partial F(x)}{\partial x_\xi} (u) \left(\psi (u)-\int\limits_\KK \psi(v)x_\xi (dv) \right) x_\xi (du).
\end{equation}

For migration transition probability $a(\cdot,\cdot)$ and migration rate $c$
\begin{equation}\label{AG-e249}
G^{\rm mig}(F)=c \sum\limits_{\xi \in \GG} \sum\limits_{\xi^\prime \in \GG} a(\xi,\xi^\prime) \int\limits_\KK  \frac{\partial F(x)}{\partial x_\xi} (u) (x_{\xi^\prime}-x_\xi)(du).
\end{equation}

All the $(G,\cF,\delta_x)$-\emph{martingale problems} are for $x \in E$ \emph{wellposed} and specify the measure-valued Markov processes
\begin{equation}\label{e332}
X=(X(t))_{t \geq 0} \; , \; E=\cM(\GG \times \KK)
\end{equation}
for every initial law on $E$.
The most basic system are interacting Fleming-Viot diffusions, i.e. $G=G^{\rm mig} + G^{\rm res}$.

\section[Description of genealogies]{Description of genealogies}\label{AG-s.descgen}

The first next and new point is now to construct for the mechanisms from the paragraph "models", the corresponding \emph{genealogical} processes as Markov process on some Polish state space. A survey of our approach to genealogies is in \cite{DG18evolution}.

\subsection[The state space of genealogies]{The state space of genealogies}\label{AG-ss.state}
The first point is to specify a \emph{state space} which allows to describe the state of the genealogy of the current population derived from the family tree respectively marked family tree grown up to the present time.
This state space should be a Polish space so that the theory of stochastic processes can be used.

\begin{figure}[t]
\begin{center}
      \includegraphics[width=0.5\textwidth]{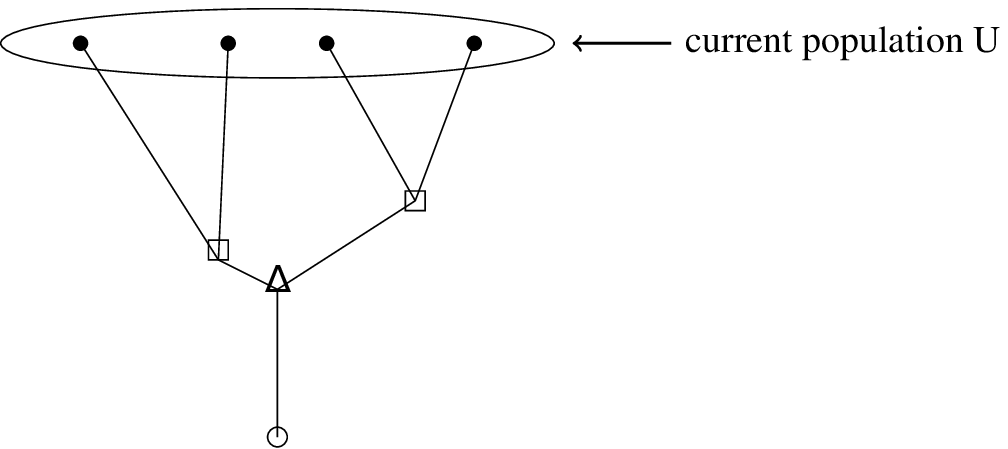}
\end{center}
\caption{\label{AG-F1}}
\end{figure}

The idea is to describe the genealogy of the population alive at time $t$ by \emph{ultrametric probability spaces}.
Fix a time $t$ referred to as the presence. Suppose we now have a family tree of the form: \textbf{Figure~\ref{AG-F1}}
\medskip

\medskip

Then define an ultrametric by setting for two \emph{currently alive individuals}: the \emph{genealogical distance} is twice the time we have to go back in time to find the most recent common ancestor $\square$ resp. $\Delta$.
Label the current individuals by $U$ and denote the genealogical distance on $U \times U$ by $r$, to get the \emph{ultrametric space} $(U,r)$.

Suppose the number of current individuals is huge, then we can study the object $(U,r)$, an complete separable ultrametric space only by drawing \emph{samples}, this is in particular the case if we pass to \emph{infinite} population size.
For taking samples we need a probability measure $\mu$ on $(U,\cB(U))$, where $\cB$ is the Borel field on $U$ w.r.t. the topology induced by the metric $r$.
We end up as objects for our description of the genealogy with \emph{ultrametric probability spaces $(U,r,\mu)$}.

Since we are not interested in the individuals as such but rather their relations, we pass to \emph{equivalence classes} w.r.t. maps $\varphi:U \to \wt U$ with $\wt r(\varphi(u),\varphi(u^\prime))=r(u,u^\prime)$ for $\mu \otimes \mu$-almost all $(u,u^\prime) \in U \times U$ and $\varphi_\ast \mu = \wt \mu$, i.e. \emph{measure preserving isometries} between $(U,r,\mu)$ and $(\wt U,\wt r,\wt \mu)$. 
The equivalence classes are denoted:
\begin{equation}\label{AG-e245}
\cU=[U,r,\mu].
\end{equation}

The set of all such equivalence classes we denote by
\begin{equation}\label{AG-e251}
\UU_1.
\end{equation}
We turn $\UU_1$ into a \emph{Polish space} by introducing the following \emph{topology}, which arises from the idea that $\cU_n=[U_n,r_n,\mu_n]$ converges to $[U,r,\mu]$, if all sampled \emph{finite} subfamilies $\{u_1,u_2,\cdots,u_m\}$, with $m \in \NN$ from the $\mathfrak{U}_n$, which is completely described by $(r(u_i,u_j))_{1 \leq i<j \leq m} \in \RR^{m \choose 2}$, converge in law w.r.t. the push forward of $\mu^{\otimes m}$ under $r$, to the corresponding object of $[U,r,\mu]$.
Formally introduce test functions called \emph{polynomials} on $\UU_1$:
\begin{equation}\label{AG-e256}
\Phi^{n,\varphi} \left([U,r,\mu]\right)=\int\limits_{U^n} \varphi \left((r(u_i,u_j))_{1 \leq i<j\leq n}\right) \mu^{\otimes n}(d \underline{u}),
\end{equation}
where $\underline{u}=(u_1,\cdots,u_n)\in U^n,\varphi \in C_b(\RR^{n \choose 2},\RR), n \in \NN$.\\

Then abbreviating $\cU=[U,r,\mu]$, we say (here $\cU_n,\cU \in \UU_1$),$\cU_n \to \cU$ for $n \to \infty$, iff
\begin{equation}\label{AG-e261}
\Phi^{m,\varphi} (\cU_n) \to \Phi^{m,\varphi} (\cU)  \quad \mbox{for} \quad
\varphi \in C_b(\RR^{m \choose 2},\RR), \quad \forall \; m \in \NN.
\end{equation}
It has been proven in \cite{GPW09} that this topology is \emph{Polish} and a metric has been exhibited.
This topology is called the \emph{Gromov-weak topology}.
Therefore we can use the Polish space $\UU_1$ as the state space for random genealogies denoted $\mathfrak{U}$.

We can introduce now \emph{marks}, where we fix a complete separable metric mark space $(V,r_V)$ to code location, type and internal state of individuals.
This means we replace the basic set $U$ by $U \times V$, where here $V=\GG \times \KK \times \II$, with $\cB(U \times V)$ w.r.t. $r \otimes r_V$.
This gives us marked ultrametric measure spaces $(U \times V,r,\nu)$.
The measure $\mu$ from above is replaced by $\nu$ on $\cB(U \times V)$. 
The equivalence class is w.r.t. measure and mark preserving isometries of $(U \times V,r,\nu)$ and $(U^\prime \times V,r^\prime,\nu^\prime)$, i.e. with $\nu=\mu \otimes \kappa, \nu^\prime =\mu^\prime \otimes \kappa^\prime$, $\kappa$ resp. $\kappa^\prime$ transition kernel from $U$ resp. $U^\prime$ to $V$ we have $\varphi_\ast \mu=\mu^\prime, \kappa^\prime(\varphi(u),\cdot)=\kappa(u,\cdot) \; \mu-$ a.s. for $\varphi$ an isometry on $\wt U \times V \to \wt U^\prime \times V$ with $\wt U,\wt U^\prime$ the supports of $\mu$ respectively $\mu^\prime$. 
The space of all equivalence classes $[U \times V,r,\nu]$ is denoted
\begin{equation}\label{AG-e286}
\UU^V_1.
\end{equation}
The topology is defined by defining \emph{convergence of sequences} and is again given in terms of \emph{polynomials}, but where in \eqref{AG-e256} we now have in the integral of $\varphi g$ w.r.t. $\nu^{\otimes n}$ with $g \in C_b(V^n,\RR)$ hence the measure $\mu$ on $U$ is replaced by $\nu$ on $U \times V$.

The concept can be extended to finite measures, compare \cite{DG18evolution}, considering
\begin{equation}\label{AG-e358}
(\bar \cU,\widehat \cU),\bar \cU \in [0,\infty), \widehat \cU \in \UU_1,
\end{equation}
taking the product topology and identifying here all $(0,\cU)$ to get the spaces
\begin{equation}\label{AG-e353}
\UU, \; \mbox{respectively}\; V\;\mbox{-marked}\; \UU^V.
\end{equation}

The concept has to be extended further to cover \emph{infinite} measures on $U \times V$ which are only restricted to $\GG_n \times \KK \times \II$ finite where $\GG_n \uparrow \GG$ and $\GG_n$ are \emph{finite} or \emph{bounded} sets.
In this case we consider "the sequence of all the restrictions to $\GG_n$" as the basic objects which converges if every element of the sequence does (see \cite{GSW} for details).

\subsection[The dynamics of genealogies]{The dynamics of genealogies}\label{AG-ss.dynamic}
The populations and their genealogy we study follow a stochastic process $\mfU=(\mfU_t)_{t \geq 0}$, which is characterized by a \emph{martingale problem}, which is specified by an \emph{operator} on a domain of test functions the \emph{polynomials}.
We specify now the actions of the operators on polynomials, the operators corresponding to the listed mechanisms (below \eqref{AG-e179}) but they are now driving an $\UU^V$-valued process.
They are calculated by considering the \emph{individual based process}, a simple jump process with deterministic growth of distances, where the generator is easily calculated and then we take the many individuals, small mass (and sometimes scaled rate) limit.

We consider as test functions now polynomials $\Phi=\Phi^{n,\varphi,g}$ of degree $n \geq 1$.
We obtain the operator $\Omega$ of the \emph{$\UU^V_1$-valued Fleming-Viot process} $\mathfrak{U}$ with selection and mutation now as follows:
\begin{equation}\label{AG-e357}
\Omega=\Omega^{\rm grow} + d \Omega^{\rm res} + s \Omega^{\rm sel} + c \Omega^{\rm mig}.
\end{equation}
(1) Then the basic single terms describing the evolution of the \emph{genealogy}  are now:
\begin{equation}\label{AG-e363}
\Omega^{\rm grow} (\Phi^{n,\varphi,g})=\Phi^{n,2 \bigtriangledown  \varphi,g} ,\quad \bigtriangledown \varphi(\uur)=\sum\limits_{1 \leq i<j\leq n} \frac{\partial \varphi(\uur)}{\partial r_{i,j}},
\end{equation}

\begin{equation}\label{AG-e369}
\Omega^{\rm res} (\Phi^{n,\varphi,g})=\frac{1}{2} d \; \Phi^{n,\wt \varphi,g}, \quad \wt \varphi= \sum\limits_{1 \leq i<j\leq n} (\varphi \circ \theta_{i,j}-\varphi) (g \circ \wt \theta_{i,j}),
\end{equation}
where the replacement operators $\theta_{i,j},\wt \theta_{i,j}$ are  defined as:
\begin{equation}\label{AG-e373}
\theta_{i,j} (\uur)=\wt \uur, \; \wt r_{k,\ell}= \left \{ 
\begin{array}{lcl}
r_{k,\ell} & if & k,\ell \neq \{i,j\} \\
r_{i \wedge k, i \vee k} & if & j=\ell \\
r_{j \wedge k,j \vee k} & if & i=\ell,  
\end{array}\right.
\end{equation}
\begin{equation}
\wt \theta_{i,j} (\underline{v})=\wt v,\; \wt v_k = \left \{ 
\begin{array}{lcl}
v_k & if & k \notin \{i,j\} \\
v_i & if & k \in \{i,j\}.   
\end{array}\right.
\end{equation}
The operator $\Omega^{\rm grow}$ describes the \emph{aging}, which means that distances grow with time at rate $2$.
The operator $\Omega^{\rm res}$ describes the \emph{splitting} in the tree where two new branches marked with their type and location start growing.
We can add here also $\Omega^{\rm Cann}$, see \cite{GKW}, but we do not write this out here.
\medskip

\noindent (2) The other operators describe the evolution of the \emph{marks} on the genealogy, $\Omega^{\rm mig}$ the change of location and the $\Omega^{\rm mut},\Omega^{\rm sel}$ the changes of types which occur and in the latter case the insertion of an individual into the sample.

Assume here that $V=\GG$ for convenience. 
The operator $\Omega^{\rm mig}$ is defined on $\Phi^{n,\varphi,g}$ as follows
\begin{equation}\label{AG-e403}
\Omega^{\rm mig} \Phi^{n,\varphi,g} =\sum\limits^n_{k=1} \Omega_k^{\rm mig} \Phi^{n,\varphi,g}_{k,\xi,\xi^\prime},
\end{equation}
where
\begin{equation}\label{AG-e407}
\Omega^{\rm mig}_k \Phi^{n,\varphi,g} = \sum\limits_{\xi,\xi^\prime \in \GG} a(\xi,\xi^\prime) \left(\Phi^{n,\varphi,g_k^{\xi,\xi^\prime}}-\Phi^{n,\varphi,g}\right)
\end{equation}
and $\Phi^{n,\varphi,g}_{k,\xi,\xi^\prime}=\Phi^{n,\varphi,g_k^{\xi,\xi^\prime}}$  with
\begin{equation}\label{AG-e411}
g_k^{\xi,\xi^\prime} (v_1,\cdots,v_{k-1},\xi,v_{k+1},\cdots,v_n) =g (v_1,\cdots,v_{k-1},\xi^\prime,v_{k+1},\cdots,v_n).
\end{equation}
In the general case where $g$ depends on a mark with a location and type component then $\Omega^{\rm mig}$ acts only on the $\GG$-component of the mark.

The mutation operator $\Omega^{\rm mut}$ work similar as $\Omega^{\rm mig}$, it just acts on a function $g$ of the type and $a(\cdot,\cdot)$ is replaced by the mutation kernel $m(\cdot,\cdot)$.

Finally we come to $\Omega^{\rm sel}$. Set (selection based on a fitness potential $\psi$):
\begin{equation}\label{AG-e417}
\Omega^{\rm sel} \Phi^{n,\varphi,g} = \Phi^{n,\varphi,\wt g}, \mbox{  with  } \wt g(\underline{u},u_{n+1})=\sum\limits^n_{k=1} g(\underline{u})\psi(u_k)-g(\underline{u})\psi(u_{n+1}).
\end{equation}
Note that here the \emph{degree} of the monomial \emph{increases by one}, this represents the additional potentially fitter type being added to the sample, this effect we saw for the measure valued version already, this means we build up a decision tree to see whether or not a fitter individual replaces the $k$-individual or not.

In \cite{DGP12,DGP13,GSW,GPW09} we showed that the various combinations of these operators all lead to \emph{wellposed martingale problems} on $\UU^V$, the process is denoted
\begin{equation}\label{e480}
(\mfU_t)_{t \geq 0}.
\end{equation}

This setup allows us to model also the \emph{genealogies of sweeps} in spatial systems with \emph{rare mutation and selection} in a fashion,  which was laid out in the monograph \cite{DGsel14} on the measure-valued process, to analyze the configurations of types as they evolve through rare fit mutants to subsequently fitter populations.
We cannot review this large amount of material here, but note that it remained open to study the \emph{genealogies} corresponding to the effects found  there.
The foundations needed for doing this have now been developed in \cite{DGP12,DGP13} and \cite{GSW} and have now to be used in forthcoming work.

The reader might have noticed, we have not discussed \emph{recombination} at this point.
Here however new foundational work is needed as we need \emph{multi-metric} measure spaces to cope with several lines of descent from one individual. 
This of course requires new type of martingale problems, dualities and test functions.
The investigations here (Depperschmidt, Greven, Pfaffelhuber, see \cite{DGP19}) are not yet complete.

\section[HMFL, continuum space limit and finite system scheme]{Hierarchical mean-field limit (HMFL), continuum space limit and finite system scheme}\label{AG-s.hierlim}

Here we describe some principal methods to study effects in spatial systems and to investigate \emph{universal large-scale properties} of model classes. These methods are as the \emph{hierarchical mean-field limit} or the \emph{continuum space limit} and the \emph{finite system scheme}.

The big questions which these tools allow to answer is:
\begin{itemize}
\item What are the \emph{universality classes} for the large scale behaviour in spatial population models?
\item What does the theory of populations in \emph{infinite} geographic space have to do with that of systems on \emph{large finite} spaces, we may see in reality?
\end{itemize}

The first point is tackled using two different approaches.
For one the \emph{multi-scale renormalization} combined with the hierarchical mean field limit and second the continuum space limit.
The second point is approached via the finite system scheme.

\subsection[Renormalization via hierarchical mean-field limit (HMFL)]{The hierarchical mean-field limit (HMFL)}\label{AG-ss.basic}

We introduce here a specific form of large space-time scale \emph{renormalization analysis} based on the \emph{hierarchical mean-field limit}.

\paragraph{The basic idea of the hierarchical mean-field limit} 
Here we explain the idea to replace, for our spatial systems $X$ (measure-valued) or the $\UU^V$-valued genealogy processes $\mfU$, the geographic space $\ZZ^2$ by a different abelian group with \emph{similar potential theoretic properties} but simpler structure.
This leads us to a geographic space which had been suggested for ecological models for different reasons namely to be more appropriate then an euclidean geographic structure.
But it is in fact also more simple to treat due to the tree structure, namely the \emph{hierarchical group}:
\begin{align}\label{AG-e275}
\Omega_N= \bigoplus\limits_{i \in \Nnull} Z^i_N, \quad Z^i_N \; \mbox{is the cyclical group} \; \{0,1,\cdots,N-1\},
\end{align}
with $+$ the addition modulo $(N)$.
Elements of $\Omega_N$ are denoted $\underline{i}$.

We can view $\Omega_N$ as the leaves of an \emph{$N$-ary tree} which suggests to use the following ultrametric on $\Omega_N$:
\begin{align}\label{AG-286}
d(\underline{i},\underline{j}) & = k \; \mbox{for} \; \underline{i}=(i_0,i_1,\cdots),\underline{j} =(j_0,j_1,\cdots) \in \Omega_N \; \mbox{iff} & \\
i_\ell & = j_\ell \; \mbox{for} \; \ell \geq k \; \mbox{and}\; i_{k-1} \neq j_{k-1}, \mbox{  or  } 0 \mbox{  for  } \underline{i}=\underline{j}.& \nonumber
\end{align}

\begin{figure}[t]
\begin{center}
      \includegraphics[width=0.5\textwidth]{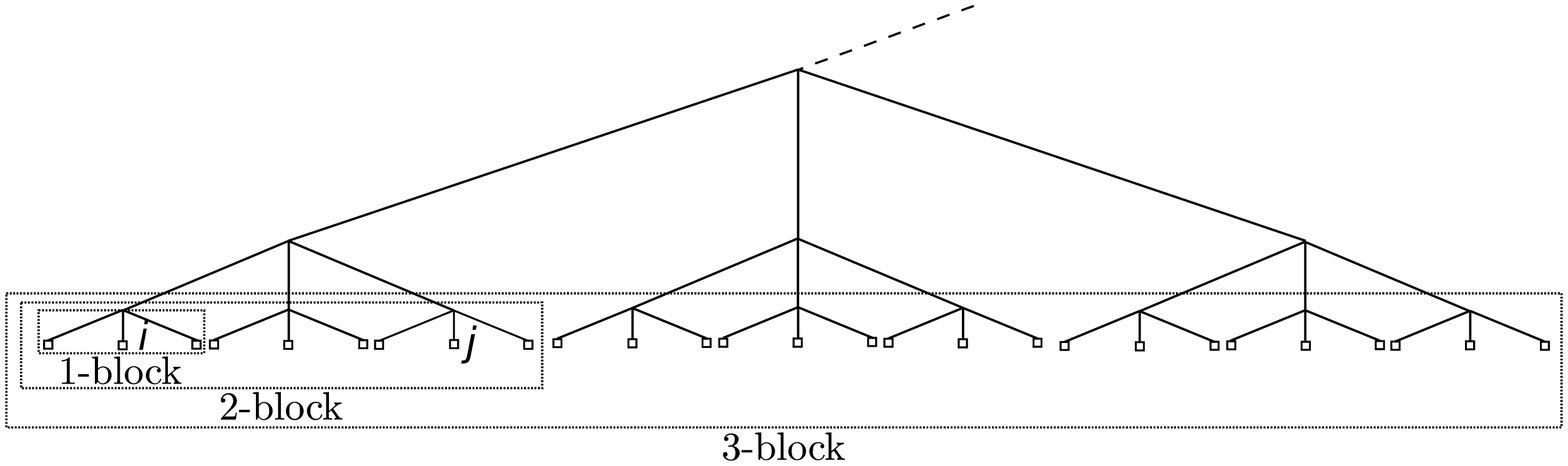}
\end{center}
\caption{\label{AG-F2} }
\end{figure}

\textbf{Figure~\ref{AG-F2}}

\medskip

This object has for another ultrametric $\wt d$ which is measuring the distance by the volume of the ball with this radius, i.e. $\wt d(\underline{i},\underline{j})=exp(d(\underline{i},\underline{j}))$, the Hausdorff dimension $log N$, which is of course not $2$.
However the potential theoretic properties of a geographic space depend on its structure \emph{and} the structure of migration and this interplay can be described by the \emph{degree} of the migration on a group.

We consider here $a_N(\ui,\uj)=a_N(0,\uj-\ui)$ arising from a sequence $(c_k)_{k \in \Nnull}$ as follows.
We pick a ball with diameter $\wt d=e^k$
\begin{align}\label{AG-e483}
&\mbox{ at rate $c_k/N^k$}  
\end{align}
and then jump to a point according to the uniform distribution on that ball.

The potential properties of the pair, $\Omega_N$ and this migration,  are best described by the \emph{degree} $\gamma \in [-1,\infty)$ which corresponds to strong recurrence for $\gamma \in [-1,0)$, null recurrence  for $\gamma=0$, transience for $\gamma > 0$ and strong transience for $\gamma >1$.

The degree $\gamma$ is defined by
\begin{equation}\label{AG-e298}
\gamma = sup \left\lbrace \zeta > -1, G^{\zeta +1} f < \infty \; \mbox{ for all $f$ bounded with compact support.}\right\rbrace,
\end{equation}
where 
\begin{equation}\label{AG-e302}
G f=\frac{1}{\Gamma (\zeta)} \cdot \int\limits_0^\infty t^{\zeta-1} a_t \; f \; dt.
\end{equation}

The \emph{degree} of the random walk on the hierarchical group $\Omega_N$ with $c_k=c^k$ is given by
\begin{equation}\label{AG-e307}
\gamma_{N,c}=\frac{log c}{log N-log c}.
\end{equation}
Therefore we see that $\gamma_{N,c} \uparrow 0, \gamma_{N,c} \downarrow 0$ for $N \to \infty$ if $c < 1$ resp. $c > 1$ and $\gamma_{N,c}=0$ for $c=1$ for all $N \in \NN$.
Hence we see that we can approximate for $N \to \infty$ the degree $=0$ from above and below.

Simple random walk on $\ZZ^d$ has degree $\gamma_d=\frac{d}{2}-1$, i.e. in $\ZZ^2$ the degree $0$.
Therefore we expect to be able to \emph{explore the situation on the euclidean space and in particular $\ZZ^2$ using the hierarchical group and letting $N \to \infty$}.

\paragraph{The HMFL and the interaction chain}\label{AG-HMFL}

\textbf{(1) Renormalizing measure valued processes} \quad
Consider first our \emph{measure-valued} systems $X$.
In order to carry out here the HMFL we have to give a \emph{sequence} of \emph{time-space scales} with which we \emph{renormalize} a spatial system to obtain in the limit as characterizing element the so called \emph{interaction chain}.

For that purpose we consider the \emph{averages} of the system in a sequence of $k$-balls $B_k$ (centered at $0$ say) for $k=0,1,\cdots,j+1$ and then we \emph{speed up time} by $\beta_k(N)$, here $\beta_k(N)$ is the ball volume $N^k$.
This results for fixed $t >0$ in a sequence of renormalized systems of $X=(x_{\underline{i}})_{\underline{i} \in \Omega_N}$:
\begin{equation}\label{AG-e439}
\left(\wt M_k^{N,j}\right)_{k=j+1,j,\cdots,0} : \wt M^{N,j}_k= \frac{1}{N^k} \sum\limits_{\underline{i} \in B_k} x^N_{\underline{i}} (t\beta_k(N)).
\end{equation}

The scaling limit of the renormalized systems as $N \to \infty$ is the 
\begin{equation}\label{AG-e444}
\mbox{\emph{interaction chain}} \; \left(\wt M^{\infty,j}_k\right)_{k=j+1,\cdots,0}.
\end{equation}

The basic theorem needed then is the \emph{existence} of the limit and the next step is the identification as a certain Markov chain and the  analysis of its properties.
The existence and identification has been established in many systems,  for example, interacting Fleming-Viot with selection and mutation \cite{DGsel14}, the Cannings processes with blockresampling \cite{GdHKK} and others.

Properties of the system as the \emph{dichotomy} coexistence versus formation of monotype spatial clusters have been shown in neutral Cannings or Fleming-Viot processes, establishing the dichotomy by showing that depending on the parameters of migration $\cL[M_\cdot^{\infty,j}]$ converges as $j \to \infty$ to a nontrivial \emph{entrance law} of the Markov chain corresponding to coexistence of types, respectively a trivial one i.e. monotype states. 
Here it is convenient to let the running index in \eqref{AG-e439} be in the negative integers running from $-(j+1)$ up to $0$, i.e. set
\begin{equation}\label{AG-e585}
M^{\infty,j}_k=\wt M^{\infty,j}_{-k} \quad , \; k=-(j+1),\cdots,0.
\end{equation}
These results we shall describe in Section~\ref{AG-ss.results}.
\medskip

\textbf{(2) Renormalizing $\UU^V$-valued states}\quad 
In order to study also the evolution of \emph{genealogies}, we have to consider the aggregated genealogies associated with the individuals marked with a $k$-block of a whole hierarchy of such blocks, zooming in on one point.
We focus on the state 
\begin{equation}\label{e624}
\mathfrak{M}^N_k=\left[U \times B_k(0),r_k,\mu |_{U \times B_k} \right]
\end{equation}
the macroscopic genealogy on level-$k$.
Then the process is viewed in a sequence of $k$-blocks in a time-rescaling which then requires an adapted rescaling of mass as before and then also a \emph{rescaling of distances} which now should be measured based on macroscopic time.

This leads then to a chain $(\widetilde{\mathfrak{M}}^{N,j}_k)_{k=j+1,1,\cdots,0}$ of such blocks and then letting $N \to \infty$ to a chain $(\widetilde{\mathfrak{M}}^{\infty,j}_k)_{k=j+1,1,\cdots,0}$.
The first point is of course to show the existence of this limit and then to give a description of the chain in $k$ with $k$ running from $(j+1)$ to $0$.
Next we want to identify this object as a Markov chain and finally identify its transition law.
Again as in \eqref{AG-e585} we switch and are using now a new index and get $k \to -k$
\begin{equation}\label{AG-e613}
(\mathfrak{M}_k^{\infty,j})_k=-(j+1),\cdots,0.
\end{equation}
Then we are ready to investigate existence and properties of the $j \to \infty$ limit, giving insight in the equilibrium genealogy.

\subsection[Continuum-space limit]{Continuum-space limit}\label{AG-ss.conlim}

Another method to determine universality classes of the large space-time scale behaviour is to embed the systems on $\ZZ^d$ or $\Omega_N$ in $\RR^d$ respectively $\Omega^\infty_N= \bigoplus\limits_\ZZ Z_N$ and then scale time and space such that the underlying random walks converge to a limit process on "continuum" space $\RR^d$ or $\Omega^\infty_N$ and then with possibly a further rate scaling try to obtain as the scaling limit a specific process on continuum spaces for a large class of discrete space models.
The most prominent example is the Dawson-Watanabe or super process on $\RR^d$, which appears as \emph{universal scaling limit} for many types of spatial models on $\ZZ^d$ or other discrete graphs.
This is developed in great detail for two-type branching models with family interaction in \cite{DGZ}.
We come to genealogy valued i.e. $\UU^\GG$-valued version in Section~\ref{AG-sss.contlim}.

\subsection[The finite system scheme]{The finite system scheme}\label{AG-ss.finsyst}

\paragraph{The finite system scheme $\UU^\GG$-valued}
The finite system scheme considers in spatial systems for the geographic space $\GG$ a sequence of groups $\GG_n \subseteq \GG$ with $|\GG_n| < \infty$ and $\GG_n \uparrow \GG$ as $n \to \infty$.
For example in $\ZZ^d$ we take $\GG_n=\ZZ^d \cap [-n,n]^d$ for $n \in \NN$ and to turn this into an abelian group, consider the addition modulo $2n$. 
Then choose as migration the original walk, but observe here the  position modulo $(2n)$.
These system then approximate the system on $\GG$ in the sense that if we start $X^n$ on $\GG_n$ in $X|_{\GG_n}$, then
\begin{equation}\label{AG-e200}
\cL \left[\left(X^n(t)\right)_{t \geq 0} \right] \nto 
\cL \left[\left(X(t)\right)_{t \geq 0} \right].
\end{equation}
A similar relation holds for $\UU^\GG$-valued processes, namely 
\begin{equation}\label{AG-e602}
\cL\left[(\mfU^n(t))_{t \geq 0}\right]\Longrightarrow \cL \left[(\mfU(t))_{t \geq 0}\right],
\end{equation}
where the state $\mfU^n$ arises by projecting the state $\mfU$ via projection $U \times \GG \to U \times \GG_n$.
Here we could also consider types in a set $\KK$ then of some we take $U \times (\GG \times \KK) \to U \times (\GG_n \times \KK)$.

Consider now on measure-valued systems first.
Then we observe the finite system globally respectively locally via ($\pi_W$ is the projection on a local window $W \subseteq \GG, |W| < \infty$):
\begin{equation}\label{AG-e206}
\cE_n(t)=\frac{1}{|\GG_n|} \sum\limits_{\underline{i} \in \GG_n} \delta_{x_{\underline{i}} (t)} \; \mbox{and  } \qquad \pi_W \mu_n(t) \; \mbox{with}\; \mu_n(t)=\cL \left[X^n(t)\right].
\end{equation}
Next look for a \emph{scale $\beta(n)$} such that:
\begin{equation}\label{AG-e210}
\cL \left[\left(\cE_n \left(t \beta(n)\right)\right)_{t \geq 0} \right] \nto \cL \left[\left(\mu(t)\right)_{t \geq 0}\right] \; , \; \cL \left[ \left(\mu_n\left(t \beta(n)\right)\right)_{t \geq 0} \right] \nto \cL \left[\left(\mu(t)\right)_{t \geq 0}\right],
\end{equation}
and the $\mu(t)$ are \emph{random extremal invariant measures} of the \emph{infinite} system.
More precisely, suppose we can find a statistic $\widehat \theta_n$ of the finite system such that as $n \to \infty$:
\begin{equation}\label{AG-e215}
\cL \left[\widehat \theta_n\left(t\beta(n)\right)_{t \geq 0}\right] \Longrightarrow \cL \left[(\Theta(t))_{t \geq 0}\right],
\end{equation}
where $\Theta$ is such that $\mu(t)=\nu_{\Theta(t)}$ is an extremal invariant measure of our process and the $\Theta$-process is a \emph{Markov} process.
If \eqref{AG-e206}-\eqref{AG-e215} all hold, then we say the \emph{finite system scheme holds}.

\begin{example}\label{AG-e.557}
Indeed for systems $X(t)=(x_\xi(t))_{\xi \in \ZZ^d}$ of interacting Fisher-Wright diffusions on $\ZZ^d$, with resampling constant $d$ and  with migration kernel $a(\cdot,\cdot)$ such that $\widehat a$ is transient $(\widehat a(\xi,\xi^\prime)=1/2(a(\xi,\xi^\prime)+a(\xi^\prime,\xi))$ and an i.i.d. initial distribution on $\GG$ restricted to the set $\GG_n$, with $\GG_n=[-n,n]^d \cap \ZZ^d$ the \emph{finite system scheme holds} with:
\begin{equation}\label{AG-e221}
\beta(n)=|\GG_n|=(2n+1)^d \quad,\quad \widehat \theta^n = \frac{1}{|\GG_n|} \sum\limits_{\xi \in \GG_n} x_\xi,
\end{equation}
\begin{equation}\label{AG-e2244}
\left(\Theta(t)\right)_{t \geq 0} \; \mbox{is given via} \; d\Theta(t)=\sqrt{d^\ast \Theta (t)(1-\Theta(t))} \; dw(t),\Theta(0)=E[x_0(0)].
\end{equation}
Namely \eqref{AG-e200}-\eqref{AG-e215} hold with $d^\ast=d(1+d \wh A(0,0))^{-1}$, where $\wh A$ is the Green function of $\wh a$.
\end{example}

This is a special case of the formula that if $g$ is the original diffusion function in a component then the one of $\Theta$ is $\mathcal{F} g$ where $\mathcal{F} (g)(\theta)=E_{\nu_\theta}[g]$ with $\nu_\theta$ being the unique spatially ergodic equilibrium with $E[x_0]=\theta$, see \cite{CGSh95} for detail.
\medskip

\noindent \textbf{Mean-field finite system scheme.} \quad
We can connect the idea of the hierarchical mean-field limit and the finite system scheme, to get the \emph{mean-field finite system scheme}.
Here one replaces the $\GG_n$ by $\{0,1,\cdots,n-1\}$ and uses as random walk the kernel with the uniform distribution. As infinite system we use the \emph{McKean-Vlasov limit} of the mean-field system, i.e. the limit $N \to \infty$.
Then the \emph{mean-field finite system scheme} holds iff~\ref{AG-e200}-~\ref{AG-e215} hold with these choices.
This is the key feature needed for the hierarchical mean-field limit to work.
For the interacting Fleming-Viot diffusions the mean-field finite system scheme holds.
\medskip

\noindent \textbf{Finite system scheme for genealogies.} \quad
For $\UU^V$-valued systems describing the genealogies we modify this as follows.
The first relation in \eqref{AG-e2244} is now modified, the second remains but is now the law of the genealogy of the population with values  in $\UU^\KK$.
For the first object in \eqref{AG-e206} we consider the genealogical relations on a \emph{macroscopic scale} by considering for a state $\mfU$ the scaled state:
\begin{equation}\label{AG-e635}
\mfU^{(n)}(t)= \left[U \times \GG_n,\widetilde{r}, \nu(\cdot \times \GG_n \times \cdot)/|\GG_n|\right], \mbox{  with   } \widetilde{r} (\cdot,\cdot)=r(\cdot,\cdot)/|\GG_n|.
\end{equation}

We then say that the finite system scheme holds if
\begin{equation}\label{AG-e657}
\cL \left[\mfU^{(n)} (t)\right] \nto \cL \left[\mfU^\ast(t)\right] \; , \; \cL \left[\mu_n(t\beta(n))\right] \Rightarrow \cL\left[(\mu(t))_{t \geq 0}\right],
\end{equation}
where $\mfU^\ast$ is a $\UU_1^\KK$-valued process, the macroscopic evolving type-marked genealogy and where $\mu(t)$ is $\mu_{\Theta(t)}$ with $\Theta$ as in \eqref{AG-e215}.
The latter will typically be a measure valued process with values in $\cP(\KK)$.

As before this can be carried for mean-field models leading to the \emph{$\UU$-valued mean-field finite system scheme}, which is an indispensable tool for the renormalization analysis of genealogies.

\section[Key results]{Some key results}\label{AG-s.keyres}

Here we collect some of the key results we obtained combining the approach of $\UU$-valued processes, renormalization analysis and continuum space limit: in \cite{GSW}, then \cite{GdHKK,GHK17,GKW} and finally \cite{ghs1,GdHOpr1,GdHOpr2} applying the ideas sketched in Sections~\ref{AG-s.motback}-\ref{AG-s.hierlim}.

\subsection[Some new models for evolving populations]{Some new models for evolving populations}\label{AG-ss.modev}
The basic models are the spatial Fleming-Viot models on $\ZZ^d$ or $\Omega_N$ as measure-valued process with state space
\begin{equation}\label{AG-e323}
(\cP(\KK))^V, V=\ZZ^d \; \mbox{or}\; V=\Omega_N 
\end{equation}
or their corresponding evolving genealogy processes with state space
\begin{equation}\label{AG-e327}
\UU^V_1 \; \mbox{with}\; V=\ZZ^d \; \mbox{or} \; \Omega_N 
\end{equation}
and with the mechanisms we gave in Section~\ref{AG-s.motback}.

Our research was focused on giving extensions in three \emph{new} directions,
\begin{itemize}
\item allowing \emph{continuum} space,
\item adding a new \emph{Cannings mechanism}, a mechanism with a spatial structure we call Cannings mechanism with \emph{block resampling},
\item adding a \emph{coloured seedbank}, where individuals are neither subject to migration nor resampling, but remain \emph{dormant} till they wake up to become \emph{active} again and migrate or resample.
The wake up times depend on the colour of the seedbank which allows a \emph{fat tail} for the annealed laws (over the colours) of the wake-up times of typical particles becoming active.
\end{itemize}

The key question is then to what extent this \emph{changes} the \emph{qualitative properties} of the population viewed on \emph{large space-time-scales} scales.
Here are the mechanisms in detail.
\medskip

\textbf{(i)} The spatial Fleming-Viot model on $\UU^V_1$ with $V=\ZZ^d$ or $\Omega_N$ and $V=\RR$, i.e. the \emph{continuum space}.
The discrete space case is given by a wellposed martingale problem, namely the $(G^{\rm grow} + G^{\rm res},\Pi)$-martingale problem, where $\Pi$ is the set of polynomials, (see \cite{GSW} for more), which opens up new possibilities of analysis.

These models have the remarkable property that their state of the genealogy at each fixed time can be coded by a path of measure-valued processes arising from the decomposition of the state in disjoint open $2h$-balls and collecting the vector of masses of the balls and form the path by varying $h$ in $(0,t)$. 
This is opening the possibility \emph{to conclude from theorems in configuration spaces (for suitably classes) results on the $\UU_1$-valued or $\UU^V_{\underline{1}}$-valued processes} (see \cite{Grieshammer2017,MAX_Griesshammer}, where this is worked out, as for example the validity of the finite system scheme for $\UU^V_1$-valued Fleming-Viot processes).

Here we focus on the continuum space case which we will get as scaling limit of the $V=\ZZ$ case.  
The continuum space case has a rich structure and can be analyzed using \emph{graphical representations}.
\medskip

\textbf{(ii)} The \emph{Cannings model} with \emph{block resampling
} a measure-valued process on $(\cP(\KK))^\GG$ and its genealogical version on $\UU_1^V$.
Here we have to explain only the $\Lambda$-block resampling in detail,   otherwise see Section~\ref{AG-s.motback} for the generator.

This is a model on the group $\Omega_N$, where at rate $\mu_k/N^{2k}$ the components in $B_k(\xi)(k$-balls around $\xi)$ are all taken together and undergo Cannings resampling w.r.t. sampling measures $(\Lambda_k)_{k \in \Nnull}$ and then are redistributed, namely uniformly distributed at the sites of $B_k(\xi)$.
This is seen as scaling limit of a catastrophic evolution on a faster time scale, then the on site evolution.
The fast time evolution is due to \emph{catastrophic events}.
The migration $a(\cdot,\cdot)$ is given by $(c_k)_{k \in \Nnull}$ as in \eqref{AG-e483}.

\medskip

\textbf{(iii)} The Fleming-Viot model with \emph{coloured seedbank} on $\ZZ^d$.
We introduce a model with "coloured" seedbanks, mirroring that the populations are subject to different exterior conditions evolving as a Markov process.
Consider a \emph{seedbank} which is structured and has at each geographic site $\Nnull$-indexed \emph{dormant} components.
This allows us to simulate a \emph{fat tail} wake up time in a \emph{Markovian} model.

Let $(x_\xi(t))_{\xi \in \GG}$ be the active population of type one in a spatial two type Fisher-Wright model and let $((y_\xi^m (t))_{m \in \Nnull})_{i \in \GG}$ the dormant populations of type $1$, with colour of the seedbank being $m$.
The evolution is given by the following SSDE's: $Z(t)=((x_\xi(t),(y^m_\xi(t))_{m \in \Nnull})_{\xi \in G} \in ([0,1] \times [0,1]^{\Nnull})^\GG$ satisfies for $\xi \in \GG, m \in \Nnull$:
\begin{align}\label{AG-e497}
dx_\xi(t)=\sum\limits_{\xi^\prime \in \GG} a(\xi,\xi^\prime)\left(x_{\xi^\prime}(t)-x_\xi(t)\right)dt+\sqrt{bx_\xi(t)\left(1-x_\xi(t)\right)} dw_\xi(t) && \\
 + \sum\limits_{m \in \Nnull} e_m K_m \left(y^m_\xi(t)-x_\xi(t)\right)dt, && \nonumber
\end{align}
\begin{equation}\label{AG-e500}
dy_\xi^m(t)=e_m \left(x_\xi(t)-y^m_\xi(t)\right)dt.
\end{equation}
Here $a(\cdot,\cdot)$ is a spatially homogeneous transition probability and $(e_m)_{m \in \NN},(K_m)_{m \in \NN}$ are positive numbers satisfying
\begin{equation}\label{AG-e504}
\chi=\sum\limits_{m \in \NN} e_mK_m < \infty \;,\; \varrho=\sum\limits_{m \in \NN} K_m \in (0,\infty]
\end{equation}
and $(w_i)_{i \in \GG}$ are independent standard Brownian motions.
We interpret now $K_m$ as the size of the $m$-th dormant population over the active population the latter in the continuum mass limit normed to $1$.
We can replace $x \to x(1-x)$ by a function $g$ with $g(0)=g(1)=0,\;  g(x) > 0$ for $x \in (0,1)$ and $g$ locally Lipschitz.

Our system of SSDE's has a unique weak solution with continuous path which is a Feller and strong Markov process.

A key quantity is $\varrho$ which allows to distinguish the \emph{regimes} $\varrho < \infty$ and $\varrho =+\infty$. 
In this context the \emph{typical wake up time} $\tau$ is important and  characterized by
\begin{equation}\label{AG-e513}
\PP(\tau > t)=\sum\limits_{m \in \Nnull} \frac{K_m e_m}{\chi} e^{-e_m t}.
\end{equation}
By choosing the $K_m,e_m$ suitably we can produce \emph{wake up times $\tau$ with fat tails} as we shall see below, but for the state of the population we have the \emph{Markov property due to the colour labeling of the dormant states}.

\subsection[Results on some new classes of processes]{Results on some new classes of processes}\label{AG-ss.results}

In this section we describe prominent results for the various models. 
We proceed in different subsections with scaling limits of the basic models, large-scale behaviour of Cannings models with block resampling and longtime behaviour of Fleming-Viot models with seedbank.

\subsubsection{The continuum space limit of the interacting Fleming-Viot process on $\ZZ^1$}\label{AG-sss.contlim}

Here we show that the neutral $\UU^\ZZ$-valued interacting Fleming-Viot process on $\ZZ^1$ converges for a large class of underlying motions in the \emph{diffusive scaling limit} to a process one can obtain as functional from the \emph{Brownian web}.
More precisely let $(\mathfrak{U}_t)_{t \geq 0}$ be the $\UU^\ZZ$-valued Fleming-Viot process.
We now scale space, time, genealogical distances as follows:
\begin{equation}\label{AG-e579}
x \to \frac{1}{\sqrt{\epsilon}} x \; , \; t \to \frac{1}{\epsilon} t \quad,\quad r \to \epsilon \cdot r(\cdot,\cdot), \qquad \epsilon >0,
\end{equation}
in order to obtain a $\UU^\RR$-valued process
\begin{equation}\label{AG-e583}
\left(\mathfrak{U}^\epsilon(t)\right)_{t \geq 0}.
\end{equation}
This scaling turns in the limit $\varepsilon \downarrow 0$ the underlying random walks into Brownian motion if we assume that $\sum\limits_{i \in \ZZ} a(0,i)i^2 < \infty$.
We obtain in \cite{GSW} the following characterization of the \emph{Fleming-Viot universality class on $\ZZ^1$}:

\begin{theorem}[Continuum space limit of genealogies]\label{AG-th.588}
We assume that $\mathfrak{U}_0$ satisfies that $\cL[\bar {\mathfrak{U}}_0]$ is translation invariant ergodic with mean intensity $\theta < \infty$.
Then 
\begin{equation}\label{AG-e590}
\cL \left[(\mathfrak{U}^\varepsilon_t)_{t \geq 0} \right] \mathop{\Longrightarrow}\limits_{\epsilon \downarrow 0} \cL \left[(\mathfrak{U}^\ast_t)_{t \geq 0}\right].
\end{equation}
here $\mathfrak{U}^\ast$ is a $\UU^\RR$-valued Feller and strong Markov process. \qed
\end{theorem}
We obtain $\mathfrak{U}^\ast$ in \cite{GSW} from the \emph{Brownian web} in a \emph{graphical construction}. 
The Brownian web is a random set of path which makes rigorous the idea to start in every time-space point \emph{coalescing Brownian motions}.
For the long history and details of such constructions see \cite{GSW}.
Here one constructs first a system of coalescing Brownian motions on $\RR$ with finitely many initial points.
Introducing a suitable topology on path with time-index $\RR$ one can extend this construction to coalescing Brownian motion starting in a fixed countable dense subset $D$ of $\RR$.
The next step is to pass to the closure of the arising random graph and proving it is a.s. independent of the choice of $D$ and hence  obtaining the Brownian web.

We can represent the state at time $t,[U_t,r_t,\mu_t]$ using for $U_t$ an enriched copy of $\RR$, enriched by the double points of the Brownian web, as sampling measure the Lebesgue measure and identifying the $h$-balls in terms of the dual Brownian web and as \emph{measure driving the dynamic} the counting measure on the double points, where genealogies split further making it Markovian as bimeasure-valued state.

This description allows to enrich the state space to a space $\UU^V$ with $V=D([0,\infty),\GG)$, where the marks represent the \emph{ancestral lines} of the current populations i.e. the path in space going from the present individual backward through space until its birth continuing backward with the positions of the father till his birth.
These are in the present case path in the Brownian web backward in time.
This point of view has been explored in the Moran model bei P.Seidel \cite{Seidel2015} and for branching populations namely in \cite{ggr_GeneralBranching} for the super random walk compare here  Section~\ref{AG-ss.branchmod}.

Further detailed properties of the genealogies of the spatial Moran model or the spatial Fleming-Viot model are obtained in \cite{Seidel2015} and \cite{MAX_Griesshammer} and in \cite{GKW}.
Namely this construction can be carried out for super random walk, i.e. the genealogies of the system of interacting Feller diffusions analog to the interacting Fleming-Viot diffusions.

In particular one can prove the finite system scheme for $\UU^V$-valued Fleming-Viot processes for $V=\Z^d,\Omega_N$ in the case of transient symmetrized migration.

\begin{theorem}[$U^V_{\underline{1}}$-valued Fleming-Viot process: finite system scheme, \cite{MAX_Griesshammer,GKW}]\label{th.740}
\hfill
Consider the Fleming-Viot model $(\mfU_t)_{t \geq 0}$ on $\ZZ^d$ or $\Omega_N$ with $\widehat{a} (\cdot,\cdot)$ being transient.
Then the finite system holds, namely \eqref{AG-e635} and \eqref{AG-e657}, with the time-scale, the conserved quantity and macroscopic variable as in \eqref{AG-e221} and \eqref{AG-e2244}.
As limit we get $(\mfU_t^{\rm FV})_{t \geq 0}$, the $\UU_1^t$-valued Fleming-Viot process.
\end{theorem}
This can be generalized further to the $\UU^V$-valued Cannings process (\cite{GKW}), where again $\mfU^{\rm FV}$ appears as the limit in \eqref{AG-e635}.

\subsubsection{The Cannings process (with block resampling) and renormalization analysis of its large-scale behaviour}\label{AG-sss.cannproc}

We apply the hierarchical mean-field analysis to investigate universal properties of Cannings processes both on the level of the measure-valued process as for the genealogies ($\UU^V$-valued processes).

\paragraph{Measure-valued processes}
We consider the Cannings process with blockresampling on $\Omega_N$ for given sequences $(\Lambda_k)_{k \in \Nnull}$ and hence $(\lambda_k)_{k \in \Nnull}$ and $(c_k)_{k \in \Nnull}$, which we denote
\begin{equation}\label{AG-e684}
X^{(\Omega_N)}=\left((X_\zeta (t))_{\zeta \in \Omega_N}\right)_{t \geq 0}.
\end{equation}
For each $k \in \Nnull$, we take the \emph{$k$-ball averages} and view them  at times $t N^k$:
\begin{equation}\label{AG-e666}
Y^{(\Omega_N)}_{\eta,k}(t)=N^{-k} \sum\limits_{\zeta \in B_k(\eta)}X^{(\Omega_N)}_\zeta (t) \quad,\quad \eta \in \Omega_N.
\end{equation}

We need the following quantities $(d_k)_{k \in \Nnull}$ we define recursively
\begin{equation}\label{AG-e672}
d_{k+1} = \frac{c_k(\frac{1}{2} \lambda_k+d_k)}{c_k+(\frac{1}{2} \lambda_k+d_k)} \quad,\quad k \in \Nnull.
\end{equation}

We define $Z_\theta^{c_k,d_k,\Lambda_k}$, the \emph{McKean-Vlasov process}, which is obtained observing one site in the $N \to \infty$ limit of the \emph{mean-field Cannings model}, where space is $\{1,\ldots,N\}$ and migration is according to the uniform distribution and where the initial state is i.i.d. with mean $\theta$.

\begin{theorem}[HMFL and renormalization]\label{AG-th.676}

(a) For every $k \in \NN$ uniformly in $\eta \in \Omega_\infty:$
\begin{equation}\label{AG-e679}
\cL \left[\left(Y^{(\Omega_N)}_{\eta,k}(t N^k)\right)_{t \geq 0}\right] 
\Nto \cL \left[\left(Z^{c_k,d_k,\Lambda_k}_\theta (t)\right)_{t \geq 0}\right].
\end{equation}

(b) Let $t_N \to \infty$ as $N \to \infty$ with $t_N/N \to 0$ as $N \to \infty$.
Then uniformly in $\eta \in \Omega_\infty$ and $u_k \in (0,\infty)$:
\begin{equation}\label{AG-e687}
\cL \left[\left(Y^{(\Omega_N)}_{\eta,k} (N^j t_N + N^k u_k)\right)_{k=j,\cdots,0}\right] \Nto 
\cL \left[\left(M^{(j)}_{-k}\right)_{k=j,\cdots,0}\right],
\end{equation}
\begin{equation}\label{AG-e691}
\cL \left[Y^{(\Omega_N)}_{\eta,j+1} (N^j t_N)\right] \Nto \delta_\theta,
\end{equation}
with $\theta \in \cP(\KK)$ is the single site mean measure of $X^{(\Omega_N)}(0)$. 
\end{theorem}
This gives us the limit object characterizing the limit of the renormalized system, the \emph{interaction chain}, which is a time-inhomogeneous Markov chain with $M^{(j)}_{-(j+1)}=\theta \in \cP(\KK)$ and with transition kernel from level $k+1$ to $k$:
\begin{equation}\label{AG-e697}
K_k(x,\cdot)=\nu_x^{c_k,d_k,\Lambda_k}(\cdot), \quad x \in \mathcal{P} (\KK),
\end{equation}
where the r.h.s. is the equilibrium measure of $(Z_x^{c_k,d_k,\Lambda_k}(t))_{t \geq 0}$.

The interaction chain exhibits a \emph{dichotomy}, which is shifted in the parameter space compared to classical Fleming-Viot (due to blockresampling via $(\mu_k)_{k \in \NN}$) towards more \emph{clustering}. Set
\begin{equation}\label{AG-e703}
m_k=\frac{\mu_k+d_k}{c_k} \; , \; \mu_k=\frac{1}{2} \lambda_k, \; k \in \Nnull.
\end{equation}
\begin{theorem}[Dichotomy: interaction chain]\label{AG-th.706}
\begin{equation}\label{AG-e707}
\cL\left[\left(M^{(j)}_k \right)_{k=-(j+1),\cdots,0}\right] \jto 
\cL\left[\left(M^\infty_k \right)_{k=-\infty,\cdots,0}\right],
\end{equation}
where we have:
\begin{equation}\label{AG-e713}
\mbox{Clustering:} \; \cL\left[M^{(j)}_0 \right] \jto  \cL\left[\delta_U \right] \;, \; \mbox{with} \; \cL [U]=\theta,
\end{equation}
\begin{equation}\label{AG-e716}
\mbox{Local coexistence:} \;  \sup_{\psi \in C_b(\KK,\RR)} \left[E_{\cL[M_0^{(\infty)}]} [Var(\psi)] >0\right],
\end{equation}
depending on whether
\begin{equation}\label{AG-e720}
\sum\limits_{k=0}^\infty m_k=+\infty \; \mbox{or} \; \sum\limits_{k =0}^\infty m_k < \infty.
\end{equation}
\end{theorem}

The coexistence is a relevant effect even \emph{for large finite systems}.
We have
\begin{corollary}[Mean-field finite system scheme]\label{AG-cor.630}
Consider the case $\sum m_k < \infty$. 
Then with $\beta(N) = N $, $\widehat{\theta}_N$ is the empirical mean and the process $\Theta=Z_\theta^{c_0,d_1,\Lambda_1}$, the mean-field finite system scheme holds.
\end{corollary}

\paragraph{$\UU^V$-valued processes}
If we now turn to the \emph{genealogies} we first return to the discussion above \eqref{AG-e613} and see that we have to find the transition law of the interaction chain first.
This law is given by the equilibrium state of the McKean-Vlasov process on level $k$ given by the "equilibrium" denoted $\nu_\theta^{(k)}$, i.e.
\begin{equation}\label{AG-e591}
K_k(\mfu,\cdot)=\nu^{(k)}_{\bar \mfu} (\cdot) \;,\; \mfu = [U \times V,r,\nu] \mbox{  then  } \bar \mfu (\cdot)=\nu(U \times (\cdot)).
\end{equation}
What is this "equilibrium" $\nu_\theta^{(k)}$?

For this purpose we need the $\UU_1$-\emph{valued McKean-Vlasov process}, which arises as the limit of the \emph{mean-field $\UU_1$-valued process} in the limit $N \to \infty$ if we observe \emph{one}  tagged component.
The \emph{mean-field process} has geographic space $\{0,1,\cdots,N-1\}$ and the random walk is generated by the uniform distribution on this set as transition kernel.

In the limit we obtain then as $N \to \infty$ a system with components labeled by $\NN_0$, where all components are conditioned on a global parameter independent evolutions for the $\UU_1$-valued components.
This evolution is called the \emph{McKean-Vlasov process} on $\UU_1$ or $\UU^\KK_1$ in the case with types.
We denote this process by 
\begin{equation}\label{e934}
(\mfU^{\rm MV}_t)_{t \geq 0}.
\end{equation}

We need that this process has a unique equilibrium for given global parameter $\theta  \in \cM_1(\KK)$.
Since in this dynamic the distances between individuals can diverge, since there need not be a most recent common ancestor, we use the following idea.

We transform the distances by
\begin{equation}\label{AG-e607}
r(\cdot,\cdot) \longrightarrow 1-e^{-r(\cdot,\cdot)},
\end{equation}
which produces an ultrametric with values in $[0,1)$ and the value $1$ if the untransformed distance grows to $\infty$.
This allows to obtain an equilibrium for the $\UU_1^\KK$-valued McKean-Vlasov process, which on level $k$ we call $\nu_\theta^{(k)}$ with $\theta \in \cM_1(\KK)$ the global parameter.

The next step is now to let $j \to \infty$ to obtain the entrance law
\begin{equation}\label{AG-e614}
(\mathfrak{M}^{\infty,\infty}_k)_{k \in \ZZ^-}.
\end{equation}
This object describes the structure of the \emph{equilibrium measure of genealogies} for the system in the limit $N \to \infty$.
Equilibrium measure are here for the process where the metric is transformed, to avoid distances tending to $\infty$, if we have no most recent common ancestor for two individuals.

Here we obtain the same dichotomy in terms of the parameters but now on the level of $\mathfrak{M}^{\infty,0}$, this is the dichotomy between \emph{single} versus \emph{countably many founding fathers} in the equilibrium state.

\subsubsection{The Fleming-Viot process with seedbank and its longtime behaviour}\label{AG-sss.FVwseed}

Recall here that the classical Fisher-Wright diffusion exhibits the dichotomy that if we start the process in a translation invariant shift ergodic initial law with $E[x_{(0)}]=\theta \in [0,1]$, then:
\begin{equation}\label{AG-e588}
\cL [X(t)] \mathop{\Longrightarrow}\limits_{t \to \infty} \theta \delta_{\underline{1}} + (1-\theta)\delta_{\underline{0}} \quad, \quad  \underline{b} \; \mbox{the constant state} \;  (b)_{i \in \ZZ^d}
\end{equation}
versus
\begin{equation}\label{AG-e592}
\cL [X(t)] \mathop{\Longrightarrow}\limits_{t \to \infty} \nu_\theta, \quad  \mbox{where} \; E_{\nu_\theta}[x_\xi(1-x_\xi)] >0 \; \mbox{if} \; \theta \in (0,1),
\end{equation}
depending on whether $\widehat a, \widehat a(\xi,\xi^\prime)=\dfrac{1}{2} (a(\xi,\xi^\prime)+a(\xi^\prime,\xi))$, is \emph{recurrent} or \emph{transient}.

The behaviour of the population with seedbank shows a quantitative effect in the regime where $\varrho < \infty$ in the sense that the classical dichotomy of clustering versus coexistence persists but clustering occurs slower by a constant.

More interesting is the case where $\varrho=+\infty$, where recurrent migration might lead still to coexistence thus changing the behaviour qualitatively. Assume here $a(\xi,\xi^\prime)=a(\xi^\prime,\xi), \forall \; \xi,\xi^\prime \in \GG$. We assume that:
\begin{equation}\label{AG-e683}
P(\tau \in dt) \sim Ct^{-(1+\gamma)}, \; \mbox{as  } t \to \infty,
\end{equation}
which we get if:
\begin{equation}\label{AG-a687}
K_m \sim A_m \; m^{-\alpha},e_m=B m^{-\beta} \; \mbox{as  }  m \to \infty,
\end{equation}
where $A,B \in (0,\infty); \alpha,\beta \in \RR : \alpha < 1< \alpha + \beta$, when
\begin{equation}\label{AG-a691}
\gamma=\frac{\alpha+\beta-1}{\beta} \;, \; C=\frac{A}{\beta} \; B^{1-\gamma} \; \gamma \; \Gamma(\gamma).
\end{equation}
Then we get as key quantity for $\varrho=+\infty$ resp. $< \infty$:
\begin{equation}\label{AG-e695}
I_{\widehat a,\gamma}=\int_1^\infty t^{-(1-\gamma)/\gamma} \widehat a_t(0,0) dt, \quad I_{\widehat a}=\int\limits^\infty_1 \widehat a_t(0,0)dt.
\end{equation}
The seedbank \emph{favors coexistence over the classical model} as follows:

\begin{theorem}[Longtime behaviour \cite{GdHOpr1}]\label{AG-th.698}
Assume the initial distribution $\mu_0$ is translation invariant and shift-ergodic. Denote $\cL[X(t)]=\mu_t$.
\begin{itemize}
\item[(1)] Let $\varrho< \infty,\theta=E[x_0+\sum\limits_{m=0}^\infty K_m y_{m,0}/(1+\varrho)]$. Then we have

Coexistence regime: If $I_{\widehat a}<\infty$ then
\begin{equation}\label{AG-e704}
\lim_{t \to \infty} \mu_t=\nu_\theta \;,\; \nu_\theta \; \mbox{is translation invariant and shiftergodic},
\end{equation}
\begin{equation}\label{AG-a707}
\nu_\theta \; \mbox{is invariant measure}\; , \theta = E[x_0]=E[y_{m,0}], \forall \; m \in \Nnull.
\end{equation}

Clustering regime: If $\II_{\widehat a}=+\infty$, then
\begin{equation}\label{AG-e712}
\mu_t \Longrightarrow \theta \left[\delta_{(1,1^{\Nnull})}\right]^{\otimes \GG} + (1-\theta) \left[\delta_{(0,0^{\Nnull})}\right]^{\otimes \GG}.
\end{equation}
\item[(2)] Let $\varrho=+\infty$. Then the dichotomy in (1) holds if we replace $I_{\widehat a}$ by $\II_{\widehat a,\gamma}$.
\end{itemize}
\end{theorem}
This means for $\varrho = + \infty$ a shift towards \emph{more coexistence} as follows.
First of all in the \emph{critical dimension $d=2$} for the model, i.e. migration has degree $0$, we obtain now \emph{coexistence} as well as in $d=1$ provided the seedbanks are large enough. Namely
\begin{itemize}
\item[(i)] $\gamma \in (1,\infty)$: we are in case of wake-up times with finite mean, migration dominates the behaviour
\item[(ii)] $\gamma \in [1/2,1]$: interplay between migration and seedbank the effective degree of recurrence/transience of the random walk is increased if $\gamma$ decreases giving more regimes of coexistence.
\item[(iii)] $\gamma \in (0,1/2)$: seedbank dominates resulting in coexistence i.e. $I_{\widehat a,\gamma} < \infty$ for \emph{all} $\widehat a$.
\end{itemize}

In the case of coexistence and $\varrho < \infty$ the coexistence is manifest in large finite systems, as well as the reduction of volatility compared with the classical one.

\begin{theorem}[Finite system scheme]\label{AG-th.697}
For $\varrho < \infty$ and $\mathcal{I}_{\widehat{a}} < \infty$, the finite system scheme holds on $\ZZ^d$, with $\GG_n=[-n,n]^d$, $\beta(n)=(2n+1)^d, \widehat{\theta}_n=n^{-1}(\sum\limits_{\xi \in \GG_n} x_\xi+  \sum\limits_{m =0} K_m e_m y_{m,\xi})$ and the process $\Theta$ is given by
\begin{equation}\label{AG-e699}
d \Theta (t)=\sqrt{\kappa \mathcal{F}g(\Theta(t))} dw (t), \quad \mbox{with} \; \kappa=(1 + \sum\limits^\infty_{m=0}K_m)^{-2},
\end{equation}
where $(\mathcal{F}g)(\theta)=E_{\nu_\theta}[g(x_0)].$
\end{theorem}

\subsubsection{Branching models: a short outlook}\label{AG-ss.branchmod}

It is often necessary to include the possibility of \emph{fluctuating population sizes}, which means to go from Fleming-Viot models to \emph{branching} models.

Our approach to genealogies can also be used to consider branching models, for example the Feller diffusion model $dX(t)=\sqrt{dX(t)}\; dw(t)$ or the super random walk, where the $X(t)=(x_\xi(t))_{\xi \in \GG}$ with 
\begin{equation}\label{e1029}
dx_\xi (t)=\sum a(\xi,\xi^\prime)(x_\xi(t)-x_\xi(t))dt + \sqrt{bx_\xi(t)}\; dw_\xi (t), \quad \xi \in \GG.
\end{equation}
The genealogies resp. marked genealogies are $\UU$- resp. $\UU^\GG$-valued processes.
Call the first process $(\mfU^{\rm Fell}_t)_{t \geq 0}$ the second $(\mfU^{\rm SRW})_{t \geq 0}$.
The generator is obtained by replacing in \eqref{AG-e357} $\Omega^{\rm res}$ by $\Omega^{\rm bra}$ and putting $s=0$.
The operator $\Omega^{\rm bra}$ is obtained by dropping in the formula for $\Omega^{\rm res}$ \eqref{AG-e369} the term $-\varphi$ in the formula for $\wt \varphi$.
Namely one can prove that the respective martingale problems are wellposed (Theorem 1 resp. 12 in \cite{ggr_tvF14}).

These genealogical processes are closely related to the $\UU_1$ resp. $\UU_1$-valued Fleming-Viot process.
Namely
\begin{equation}\label{e1037}
\cL\left[(\wh \mfU_t^{\rm Fell})_{t \geq 0} | (\bar \mfU_t)_{t \geq 0}\right]=\cL \left[(\mfU_t^{\rm FV})_{t \geq 0} \right] \quad,\quad \cL[\bar \mfU] - a.s.
\end{equation}
where the parameter $d=b \bar \mfU_t^{-1}$, i.e. we have a time-inhomogeneous process.
The works also for the super random walk where some technicalities arise due to $0^\prime s $ in the total mass (see \cite{ggr_tvF14}).

One then proves theorems on the longtime behaviour of the super random walk in particular identifying the limit of the laws as $t \to \infty$, where again for $\wh a$ recurrent we have locally extinction and conditioned on survival mono-ancestor genealogies and for $\wh a$ transient countably many ancestors.
As we pointed out earlier this requires to introduce a transformation of the distances onto $[0,1]$ to avoid distances $+\infty$ if individuals do not have a most recent common ancestor.
Some of these issues are discussed in \cite{DG18evolution}.

We can also consider the \emph{finite system scheme} then in the context of genealogies.
This is well known for the classical super random walk.
Here the time scale $\beta_n$ is $|\GG_n|$ again and the Feller diffusion with parameter $b$ arises as the macroscopic variable $(\Theta(t))_{t \geq 0}$ in the limit $n \to \infty$.
In the case of genealogies, i.e. $\UU^\GG$-valued super random walk, the \emph{global dynamics} $(\mfU^\ast(t))_{t \geq 0}$ is the (non-spatial) $\UU$-valued Feller diffusion $(\mfU^{\rm Fell}_t)_{t \geq 0}$ and the macroscopic variable $(\Theta (t))_{t \geq 0}$ is $(\bar \mfU_t)_{t \geq 0}$ and is the classical $\RR$-valued Feller diffusion.

For the genealogy of branching models one can introduce a theory of \emph{infinite divisibility} for genealogies, i.e. $\UU$ resp. $\UU^V$-valued stochastic processes with a corresponding \emph{Levy-Khintchine formula}, based on semigroup structure developed in \cite{infdiv}.
This leads for processes with the \emph{branching property} to a decomposition into subfamilies of any prescribed depth-$h$.
This gives in turn a representation of the current state into a concatenation of a \emph{Cox point process} on the states of single ancestor depth-$h$ subfamilies.
Furthermore we obtain an identification of the associated parameters, the Cox-measures and the laws of the single ancestor subfamilies \cite{ggr_tvF14}.

This can be applied in particular to the spatial relative of the $\UU$-valued Feller diffusion, the super random walk.
Here one can obtain the cluster representation, which helps analyzing the asymptotic behaviour, see Theorem 12 in \cite{ggr_tvF14}.

Then one can prove furthermore (Theorem 4 in \cite{ggr_tvF14}) the following representation in terms of random elements in $\UU$ of diameter $2h$ corresponding to single ancestor subfamilies, which are \emph{conditioned on the time $t-h$ population size then independently} distributed according to a law $\varrho^t_h$ which has various representations.

We have the \emph{Cox cluster
    representation} (see \cite{D93} page 45/46 for that concept),
  i.e.\ a unique decomposition into \emph{depth-$h$ mono-ancestor
    subfamilies} of the time $t$ population and $\mfU_t$ can be
  represented accordingly as concatenation over a Cox point process, where Cox measure and \emph{single ancestor subfamily law} are given explicitly as law of $\bar U_{t-h}$ respectively the normalized entrance law from the zero element of $\UU$ or via a time inhomogeneous spatial coalescent.
 
  More precisely we can
  represent the $t$-top of $\mfU_t$ as concatenation of a
  $\Pois(\bar \mfU_0)$ number of random elements chosen at
  random according to $\varrho_t^t$ (decomposition in the families
  of the founding fathers).
  Hence for given
  $\bar \mfU_{t-h}=u$ we have the representation as concatenation of
  $\Pois((bh)^{-1}u)$ distributed number of independent realizations
  of $\varrho^t_h$.
  
  We obtain a \emph{decomposition in open $2h$-balls
    grouped in open $2t$-balls}. More precisely, we obtain a
  decomposition in $M_h^t$ different open $2t$-balls each of which
  is decomposed in $M_h^t$ many open $2h$-balls $\mfU^i_k$,
  where the $\mfU^i_k,k=1,\dots, M_h^t$,
  $i=1,2,\dots,M_h^t$ are \emph{independent} of $M_h^t$ and
  the $N^{t,(i)}_h$ and i.i.d.\ $\UU(h)$-valued random variables
  distributed according to $\wh \varrho^t_h$.

  Let $(Y_i)_{i \in \NN}$ be i.i.d.\ $\Exp((t-h)b/2)$-distributed, and
  let $M_h^t$ be $\Pois(\frac{2}{bh} Y_i)$ distributed and
  independent. Then we consider the $2h$-balls in distance less than
  $2t$ and group them in the $2t$-balls we obtained above. Consider
  the set $I$ of $i$ with $M_h^t\geq 1$, their number is
  $M_h^t$.
  
The $\UU^\GG$-valued Fellers branching diffusion (super random walk), actually has one very interesting feature. Namely a further point here is, that we can consider here the genealogy where we mark individuals with their \emph{complete ancestral path}.
Take the super random walk $\mfU_t=[U_t \times \GG, r_t,\mu_t]$.
We can enrich this object to a process with mark space $V=D(\RR,G)$.
Namely with each element of $U_t$ we associate as \emph{mark} a path which is constant before time $0$ and after time $t$ and follows through $\GG$ the location of the individual at earlier times back to the father where it continues with the fathers position earlier till his birth, the grandfathers position, till the position of the ancestor at time $0$.
Then we obtain an enrichment of the $\UU^\GG$-valued super random walk to the $\UU^V$-valued process $(\mfU_t^{\rm SRW,anc})_{t \geq 0}$.
This process is a branching process (generalized branching property) with infinitely divisible marginal distributions if this holds initially, see \cite{ggr_GeneralBranching} for details.

\section[Perspectives]{Future perspectives}\label{AG-s.future}

There are a number of directions to continue the research:
\begin{itemize}
\item the \emph{genealogy process} of the model with \emph{seedbank} \cite{ghs_gen},
\item description of the \emph{genealogy} in the regime of \emph{cluster formation} in two-dimensional Fleming-Viot \cite{GLW18} by analyzing the genealogy via equivalence classes of metric \emph{bi-measure} spaces,
\item the study of the \emph{spatial continuum limit} of Cannings models and their genealogies on the hierarchical group \cite{ghs1}-\cite{ghs_renorm},
\item study sweeps in the sense of \cite{DGsel14} on the level of \emph{evolving genealogies} systematically,
\item study evolving genealogies in interacting Fleming-Viot processes with \emph{recombination} \cite{DGP19},
\item study in \emph{continuum} space the Fleming-Viot model with \emph{selection} in its \emph{evolving genealogy}, in particular in the realm of rare mutation \cite{DGsel14} using the measure-valued coding of the genealogies as well as its Brownian net or Levy net coding using \cite{GSW}.
\item develop the \emph{renormalization analysis} of the $\UU^\GG$-valued super random walk and more general \emph{state dependent $\UU^\GG$-valued branching diffusions}.
\item study genealogies in branching dynamics with competition, the analog of selection in Fleming-Viot models.
A typical example is logistic super random walk.
\end{itemize}

\paragraph{Acknowledgment:} 
We thank the collaborators of many years on the described research: D.Dawson, A.Depperschmidt, P.Gl\"ode, F.den Hollander, M.Oomen, P.Pfaffelhuber, T.Rippl, R.Sun, A.Winter.\\
\hfill\\
The collaborations were supported through the SPP1590, through the DFG Grant GR 876-16.1 and 16.2, and the Grants DFG GR 876-15 and later DFG GR 876-17.

\newpage
\bibliography{spp-report}
\addcontentsline{toc}{section}{References}
\bibliographystyle{alpha}

\end{document}